\documentclass[10pt]{smfart}
\usepackage[english]{babel}
\usepackage{times,amsfonts,amsmath,amstext,amsbsy,amssymb,
  amsopn,amsthm,upref,eucal, amscd}
\usepackage[T1]{fontenc}
\usepackage[pdftex]{graphicx}
\usepackage{tcolorbox}
\newtheorem*{main lemma}{Main Lemma}
\newtheorem{theorem}{Theorem}[section]
\newtheorem{lemma}[theorem]{Lemma}
\newtheorem{corollary}[theorem]{Corollary}
\newtheorem{proposition}[theorem]{Proposition}

\numberwithin{equation}{section}
\newtheorem*{theorem*}{Theorem}
\newtheorem*{claim}{Claim}
\newtheorem*{sublemma}{Sublemma}
\newtheorem{innercustomthm}{Theorem}

\newenvironment{customthm}[1]
 {\renewcommand\theinnercustomthm{#1}\innercustomthm}
 {\endinnercustomthm}

\theoremstyle{definition}
\newtheorem{definition}[theorem]{Definition}
\newtheorem{example}[theorem]{Example}

\newtheorem{remark}[theorem]{Remark}

\def\leq{\leqslant }
\def\geq{\geqslant}

\begin{document}
\begin{otherlanguage}{english}
\title{Non-statistical behavior via Statistical instability: Non-statistical Anosov-Katok diffeomorphisms}
\author{Amin Talebi}
\address{(amin.s.talebi@gmail.com), Department of Mathematics, Sharif University of Technology, Azadi Street, Tehran, Iran and, Institute for Research in Fundamental Sciences (IPM), Niavaran Square, Tehran, Iran }

\date{January 23, 2025}
%%%%%%%%%%%%%%%%%%%%%%%%%%%%%%%%
%%%%%%%%%%%%%%%%%%%%%%%%%%%%%%%%
\begin{abstract}
Non-statistical (or historical) dynamics refer to systems where a set of points with 
positive measure, relative to a reference probability measure (typically the Lebesgue
measure on a manifold), fails to exhibit a convergent sequence of empirical measures. There are various strategies to identify maps with non-statistical behavior. In the first part of this paper, we use a strategy related to \textit{statistical instability} to demonstrate that a Baire generic map in the space of Anosov-Katok diffeomorphisms of the annulus exhibits non-statistical behavior and even maximal oscillation. The notion of statistical instability was first formalized in its general format in \cite{talebi_non-statistical_2022}, where its close connection to the existence of non-statistical dynamics was established. In the second part, we extend this framework by introducing analogous formalization and corresponding results in the case of pointwise convergence of empirical measures. Additionally, we propose a general strategy to identify non-statistical dynamics within other families of dynamical systems.
\end{abstract}

\maketitle
%\pageheight{8.5truein}
%\pagewidth{6.5truein}

\tableofcontents

\section{Introduction}

The aim of this paper is twofold; the first one is to prove a theorem regarding the Baire genericity of non-statistical behavior in the space of dissipative Anosov-Katok diffeomorphisms of the annulus. The second one is to present a pointwise convergence version of an abstract setting developed in \cite{talebi_non-statistical_2022}, where a new notion of \textit{statistical instability} was defined, and its connection to the existence of non-statistical dynamics was explored. 
\subsection{Anosov-Katok diffeomorphisms of the annulus}
In \cite{AK70} Anosov and Katok introduced a method for obtaining Lebesgue measure-preserving ergodic maps with unexpected metric properties on manifolds that admit a $\mathbb{S}^1$ free action. They considered a class of measure-preserving maps on such a manifold that can be approximated by periodic maps (like rational rotations of the torus) and prove that the set of ergodic transformations is a Baire generic subset of this space (which is an intersection of countably many open and dense subsets). Herman and Fathi in \cite{FH77} pushed forward their method to construct minimal and uniquely ergodic maps. They also proved that these maps form a Baire generic subset of the space of maps that can be approximated by periodic ones. The state of the art of the mentioned works was to use this new method of Anosov and Katok to conclude the density of such properties. Herman could also apply the Anosov-Katok method to construct exotic invariant sets for $C^\infty$ diffeomorphisms of the Riemann sphere \cite{Her86}. Here we discuss Baire generic properties of Anosov-Katok diffeomorphisms of the annulus, from a statistical point of view.\\

Let us denote the annulus $[0,1]\times\mathbb{R}/\mathbb Z$ by $\mathbb{A}$ and for $r\in [0,\infty]$ the space of all $C^r$ orientation preserving diffeomorphisms of $\mathbb{A}$ by $\text{Diff}_{+}^r(\mathbb{A})$ endowed with the $C^r-$topology. We denote the $C^r-$closure of the set of all $C^r$ diffeomorphisms of the annulus which are $C^r-$conjugate to a rotation by $\mathcal{AK}^r$ and call it the space of $C^r$ Anosov-Katok maps. We also use $\mathcal{AK}_{vol}^r$ to denote the $C^r-$closure of the set of $C^r$ volume-preserving diffeomorphisms conjugate to a rotation with a conjugacy fixing every point of the boundary. The spaces $\mathcal{AK}^r$ and $\mathcal{AK}_{vol}^r$ endowed with the induced $C^r -$topology are Baire spaces. \\

Fayad and Katok in \cite[Theorem 3.3]{FK04} obtained a description for the statistical behavior of Baire generic area-preserving Anosov-Katok diffeomorphisms of the annulus. They show that a Baire generic map in the space $\mathcal {AK}_{vol}^\infty$ has only three ergodic measures, two one-dimensional Lebesgue measures on the boundary components, and the volume measure of the annulus. In this paper, we obtain a description of the statistical properties of Baire generic dissipative Anosov-Katok diffeomorphisms of the annulus and show that they exhibit an extreme type of non-statistical behavior.
\begin{customthm}{A}\label{Anosov-Katok example}
 A Baire generic map in the space of Anosov-Katok maps of annulus $\mathcal {AK}^r$ has exactly two ergodic invariant measures each of which is supported by a different boundary component and moreover, the map is maximally oscillating.
\end{customthm}

\begin{definition}
A map $f$ is said to have \emph{maximal oscillation} if the sequence of empirical measures
\begin{equation*}
 e_{n}^{f}(x):=\frac{1}{n} \sum_{i=0}^{n-1} \delta_{f^i (x)}. 
\end{equation*}
of almost every point accumulates to every invariant measure, going to some sub-sequence.
\end{definition}

\begin{definition}\label{def.non-stat-ess} The map $f$ is called \textit{non-statistical} if there is a set of positive measure for which the sequence of empirical measures $\{e_n^f(x)\}$does not converge.
\end{definition}
\begin{remark}
Note that an Anosov-Katok map of the annulus has at least two invariant measures that are supported on different boundary components, so a generic map in $\mathcal {AK}^r$ has the minimum possible number of ergodic measures, and at the same time, maximum possible oscillation for empirical measures of almost every point, and hence it should be non-statistical. 
\end{remark}

In \cite{HK2}, Hofbauer and Keller demonstrated the existence of maximally oscillating dynamics within the family of logistic maps. Similarly, in \cite{talebi_non-statistical_2022}, the author established analogous results for the family of rational maps on the Riemann sphere. Despite great differences in their dynamical properties—such as positive entropy in logistic and rational maps and zero entropy in Anosov-Katok maps—the phenomenon of maximally oscillating maps in all three families can be attributed to a common underlying mechanism: \textit{statistical instability}.

Statistical instability, in essence, refers to the sensitivity of the statistical behavior of a dynamical system to small perturbations, leading to significant variations in statistical properties. This concept was rigorously defined and studied in \cite{talebi_non-statistical_2022}, where its deep connection to the existence of non-statistical maps was explored. One of the contributions of this paper is the combination of the abstract framework developed in \cite{talebi_non-statistical_2022} with the Anosov-Katok method. We believe that more significantly, the paper introduces a general strategy for identifying maximally oscillating—and, more broadly, non-statistical—maps across other families of dynamical systems. 
\subsection{A general formalization of statistical instability, and its connection to the existence of non-statistical dynamics}
Several results on statistical stability have been established in the context of Axiom A diffeomorphisms, expanding maps, partially hyperbolic diffeomorphisms, and non-uniformly hyperbolic maps, as can be found in \cite{Ruell97-diffsrb}, \cite{APV16}, \cite{AV02}, \cite{DOl-Via-Yang2016}, \cite{Andersson2010}, \cite{Ander-Vas2020}, \cite{yang_2021}, \cite{Mi-Cao.stat.stab.MEMC21}, \cite{Baladi2015}, \cite{Feritas05}, \cite{Alves.Carvalho.Freitas.Henon10}, \cite{Alves-Soufi.stat.stab.Rov12}, and \cite{Araujo21}. Additionally, there are a few works that address the statistical instability of dynamical systems, such as \cite{T01} and \cite{Alves-Khan19}. In all of these studies on statistical (in)stability, the dynamics are assumed to have finitely many physical measures, and statistical stability is defined as the continuous evolution of the physical measures as the dynamics change. However, as previously mentioned,
\begin{center}
    \fbox{\parbox{0.8\textwidth}{we extend the concept of statistical (in)stability and define it in a more general sense, applicable to any dynamical system, regardless of its statistical behavior. In particular, we do not assume the system has any physical measures.}}
\end{center}
 Let us now be more precise. Let $(X,d)$ be a compact metric space, $f:X\to X$ be a continuous map, and $\mu$ be a reference probability measure on $X$. By $\mathcal{M}_1(X)$ we denote the space of probability measures on $X$ endowed with Wasserstein distance $d_w$ (see section \ref{sec.abs.setting} for definitions). Suppose $\Lambda$ is a set of continuous self mappings of $X$ which is endowed with a topology finer than $C^0$ topology (not necessarily strictly finer, it may be the $C^0$ topology itself), e.g. the set of all $C^r$ diffeomorphisms of a compact smooth manifold endowed with $C^r$ topology, for any $r\geq 0$. 
\begin{definition}\label{def.stat.uns} A map $f\in\Lambda$ is called \textit{statistically unstable} if the following condition holds:
$$\lim_{N\rightarrow \infty}\limsup_{\substack{{h,g\to f}\\ {h,g\in \Lambda}}}\int_X \sup_{N<n,m}d_w(e_n^f(x),e_m^f(x))d\mu(x)>0.$$
\end{definition}

As an illustrating example, let $\Lambda$ be the set of rotations of the circle, then $f\in\Lambda$ is statistically unstable if and only if $f$ is a rational rotation.

\begin{remark}
 In an upcoming paper, we will discuss how our definition generalizes the previous definitions of statistical stability in the case where all maps in $\Lambda$ possess finitely many physical measures with the full union of basins.
\end{remark}
 
We are ready to state our main result in this direction: 

\begin{customthm}{B}\label{thm.main.essential}
Suppose $X$ is a compact metric space endowed with a reference probability measure, and $\Lambda$ is a set of continuous self-mappings of $X$, endowed with a topology finer than $C^0$ topology. Then the non-statistical maps form a Baire generic subset of the interior of statistically unstable maps. 
\end{customthm}

Theorem \ref{thm.main.essential} provides a strategy for identifying new classes of non-statistical maps. It suggests that, for this purpose, it is sufficient to find families of dynamical systems where the set of statistically unstable maps has a non-empty interior. Non-statistical logistic maps (\cite{hofbauer1990} and \cite{HK2}), non-statistical rational maps (\cite{talebi_non-statistical_2022}), and Theorem \ref{Anosov-Katok example} serve as examples of this strategy. Additionally, there is another approach to obtaining non-statistical behavior, which involves the concept of wandering domains. This idea was first introduced by Colli and Vargas in \cite{CV01}, and has since been further developed by others, including \cite{Kriki&Soma}, \cite{Pablo23}, \cite{KYT19}, \cite{KLS16}, \cite{KS15}, and \cite{Berger-Biebler2020}. One key difference between these two strategies is that the Baire genericity of non-statistical maps is ensured by statistical instability, while there are no results yet regarding the Baire genericity of non-statistical maps via wandering domains. In \cite{Kriki&Soma}, it is shown that non-statistical maps with wandering domains are dense in the Newhouse domains within the space of surface diffeomorphisms. However, it remains an open problem whether a Baire generic map in the Newhouse domains is non-statistical. We believe this problem is challenging, and our strategy of statistical instability may be useful in addressing it. In an upcoming paper, we will explore this particular issue in more detail.\\

Now, let us explain some applications of Theorem \ref{thm.main.essential}. For a dynamical system $f: X \to X$, a point $x \in X$ is called non-statistical if the sequence of empirical measures $e_n^f(x)$ does not converge. Non-statistical points appear in many well-known dynamical systems. For example, in the shift map (see Example \ref{exmpl symbolic dynamics}) and in any other dynamics that include a shift map as a subdynamics. Some of the first results in this area is the Baire genericity of non-statistical points within the phase space for subshifts of finite type, Anosov diffeomorphisms, and, more generally, any map with the periodic specification property (see \cite{Zig70}, \cite{Zig74}, \cite{HLT21}, and \cite[Theorem 1.1.11]{GK16}). In \cite{Carvalho2021SensitivityAH}, Carvalho et al. introduce a sufficient condition for the Baire genericity of non-statistical points in the phase space of a dynamical system. Although, at first glance, Theorem \ref{thm.main.essential} seems to be applicable only for obtaining the Baire genericity of non-statistical maps in a family of dynamical systems, but indeed it can be applied in a way to give the Baire genericity of non-statistical points in the phase space of a single system, and unifies most of the previous results in this direction. To achieve this, we need to define statistical stability and instability for points in the phase space. This is done in Section \ref{Sec.gen.non-stat.phase}, where we prove:

\begin{customthm}{C}\label{theorem.generic.phase}
Let $f:X\rightarrow X$ be a continuous map of a compact metric space. If every point $x\in X$ is statistically unstable, then Baire generic points are non-statistical. 
\end{customthm} 

To see how general is this theorem, note that if the map $f$ has two different periodic points and the stable set of each of them is dense in $X$, then every point in the phase space would be statistically unstable. This would be the case in many examples of hyperbolic, or partially hyperbolic maps. As another corollary, we will show:  

\begin{corollary}\label{cor.minimal.not.erg}
Let $f$ be a minimal but not uniquely ergodic map on a compact smooth manifold $X$ (such as Furstenberg's example on the torus). Then a Baire generic point in $X$ is non-statistical.
\end{corollary}

\subsection*{Questions:}
 
This paper provides tools to study the statistical behavior of generic dynamical systems in an abstract class. When the class of dynamical systems is formed by dissipative $C^r$-diffeomorphisms of a compact manifold, this study is traditionally related to the notion of physical measures.  We recall that an invariant probability measure $\nu$  is \emph{physical} if its basin $B_\nu:=\{x\in M: e _n(x)\to \nu\}$ has positive Lebesgue measure. \\

 \textbf{Question:}
For $r \geq 2$, is it true for generic $f$ in $\mathrm{Diff}^r(M)$ that the union of the basins of the physical measures of $f$ has full Lebesgue measure in M?\\

This question has been asked by Wilkinson and Shub in \cite{WS00}, but was in the mind of several other people 
 %There are also other people who care about this question
 (see \cite{Palis00}, \cite{BV00} and \cite{PS95}).\\
  
Let us now relax the conditions on physical measures and develop some questions on the abundance of non-statistical dynamics. \\

\textbf{Question:}
 Is there any reasonable non-trivial finite-dimensional family of dynamical systems that the set of parameters for which the dynamics is non-statistical has a positive measure?\\

\textbf{Question:}
Is there any open subset of dynamics in which the non-statistical maps are generic? In Newhouse domains?\\

These questions are related to the following:\\

\textbf{Question:}[Takens' last problem,\cite{Takens08}]
Can non-statistical dynamics exist persistently within a non-trivial class of smooth dynamical systems? \\

\subsection*{Appendices}

In the appendices, alternative versions of the definitions of non-statistical dynamics
and statistical instability are developed, corresponding to different types of convergence.
Specifically, for each non-negative integer $n$, $e_n^f$ defines a map from the phase space 
to the set of probability measures:
$$e_n^f:X\rightarrow \mathcal{M}_1(X).$$

Thus, as a sequence of maps from a probability space to a metric space, various types of
convergence can be considered: pointwise, $L^1$, and convergence in law. Below, we outline
three possible definitions based on these types of convergence:

\vspace{1em}
\begin{itemize}
    \item The map $f$ is called \textit{non-statistical} if the sequence of maps $\{e_n^f\}$ is 
    not essentially convergent:
    $$\lim_{N\rightarrow \infty}\int_X \sup_{N<n,m}d_w(e_n^f(x),e_m^f(x))d\mu(x)>0.$$
    (This is equivalent to Definition \ref{def.non-stat-ess}, see Lemma \ref{lemma unif}.)
    
    \item The map $f$ is called \textit{$L^1$ non-statistical} if the sequence of 
    maps $\{e_n^f\}$ does not converge in the $L^1$ topology:
    \begin{equation*}
       \limsup_{m,n\to\infty}\int_X d_w(e_n^f(x),e_m^f(x))d\mu(x)>0.
    \end{equation*}
    
    \item The map $f$ is called \textit{non-statistical in law} if the sequence 
    of pushforward measures $\hat{e}_n^f:=(e_n^f)_*(\mu)$ does not converge as a
    sequence in $\mathcal{M}_1(\mathcal{M}_1(X))$:
    $$\#\text{acc}(\{\hat{e}_n^f\}_n)>1,$$
    where $\text{acc}(\cdot)$ denotes the set of accumulation points.  
\end{itemize}
\vspace{1em}

\begin{remark}
    Recently, Coats and Melbourne, together with the author in \cite{CMT2024}, 
    demonstrated the existence of a large class of interval maps (including intermittent maps) that
    are non-statistical when considering essential convergence. However, from the perspective of
    convergence in law, these maps exhibit the best form of statistical convergence. 
    This example highlights the necessity of having distinct versions of definitions for (non-)statistical dynamics. 
\end{remark}

Similarly, we can extend the definition of statistical instability to $L^1$ 
convergence and convergence in law, obtaining counterparts to Theorem \ref{thm.main.essential} in
these settings. These developments are detailed in the appendices. Notably, the initial versions of
these definitions and results focused on convergence in law, as introduced in \cite{talebi_non-statistical_2022}.
However, since the results in this version are used for proving Theorem \ref{Anosov-Katok example}, we have included them in the appendices for the sake of completeness. Another reason for their inclusion is that some proofs presented 
here are new and more informative than those in \cite{talebi_non-statistical_2022}.

\subsection*{Acknowledgment.} 
I would like to express my special thanks to Pierre Berger and Meysam Nassiri, for their constant support and numerous comments and discussions while doing this work. I acknowledge Bassam Fayad who gave the idea of improving our result from the Baire genericity of non-statistical Anosov-Katok maps to the Baire genericity of maximally oscillating Anosov-Katok maps. I would like to list all of my financial sources while writing this paper: The Institute for Research in Fundamental Sciences (IPM), Campus France, Iranian Ministry of Science, Research and Technology, Universit\'e Paris 13 and ERC project 818737 \textit{Emergence of wild differentiable dynamical
systems}.

%%%%%%%%%%%%%%%%%%%%%%%%%%%%%%%%
%%%%%%%%%%%%%%%%%%%%%%%%%%%%%%%%%
\section{Statistical instability and non-statistical dynamics; an abstract setting}\label{sec.abs.setting}
Let $X$ be a compact metric space endowed with a reference (Borel) probability measure $\mu $ and $\Lambda $ a subset of continuous self-mappings of $X$ endowed with a topology finer than $C^0$ topology. For instance $\Lambda$ can be a subset of $C^r$ self-mappings of a smooth manifold, endowed with $C^r$ topology and $\mu$ be a
probability measure whose density w.r.t. a Lebesgue measure is a smooth positive function.
For a compact metric space $(X,d)$, Let us denote the space of probability measures on $X$ by $\mathcal M_1(X) $. This space can be endowed with weak-$*$ topology which is metrizable, for instance with \emph{Wasserstein metric} where the distance $d_w$ between two probability measures $\nu_1$ and $\nu_2$ is defined as below:
$$d_w(\nu_1,\nu_2):=\inf_{\zeta \in \pi(\nu_1,\nu_2)}\int_{X\times X} d(x,y) d\zeta \ ,$$
where $\pi(\nu_1,\nu_2)$ is the set of all probability measures on $X\times X$ which their projections on the first coordinate is equal to $\nu_1$ and on the second coordinate is equal to $\nu_2$. The Wasserstein distance induces the weak-$*$ topology on $\mathcal M_1(X)$ and hence the compactness of $(X,d)$ implies that $(\mathcal M_1(x),d_w)$ is a compact and complete metric space. 
We should note that our results and arguments in the rest of this note hold for any other metric inducing the weak-$*$ topology on the space of probability measures.\\

For a point $x\in X$ and a map $f:X \rightarrow X$, the \emph{empirical measure} 
\begin{equation*}
 e_{n}^{f}(x):=\frac{1}{n} \sum_{i=0}^{n-1} \delta_{f^i (x)} 
\end{equation*}
describes the distribution of the orbit of the point $x$ up to the $n^\text{th}$ iteration in the phase space, which asymptotically may or may not converge.
If it converges, then by observing a finite number of iterations (but possibly large) one can predict how the orbit of $x$ behaves approximately for larger iterations, from a statistical point of view. However, it may not converge. In this case, we fix the following terminology:
\begin{definition}
 For a map $f:X\to X$ we say the orbit of a point $x$ \textit{displays non-statistical behavior}, or briefly $x$ is \textit{non-statistical} if the sequence $\{ e_{n}^f(x)\}_{n}$ is divergent. 
\end{definition}

\begin{example}\label{exmpl symbolic dynamics}
This example shows the existence of non-statistical points for well-known dynamics; the shift map $\sigma$ on $X=\{0,1\}^{\mathbb Z}$. 
 Consider a point $\omega\in X$ 
  \begin{equation*}
  \omega= 0.\underbrace{0...0}_{{n_1}}\underbrace{0101...0101}_{{n_2}}\underbrace{0...0}_{{n_3}}...
  \end{equation*}
 made by putting together consecutive blocks of zeros and blocks of zeros and ones and suppose the length of $i^{\text{th}}$ block is $n_i\,$satisfying 
 $$\lim_{i\to\infty} \frac{n_i}{n_{i+1}}=0.$$
 Then it can be checked easily that $\omega$ is a non-statistical point. 
\end{example}
One can ask how large can be the set of points for which the empirical measures does not converge. The notion of non-statistical dynamics arises here. See Definition \ref{def.non-stat-ess}.\\

\begin{remark}
    Being non-statistical for a map $f:X\rightarrow X$  really depends on the reference measure. In Example \ref{exmpl symbolic dynamics}, if the reference measure is taken to be the Bernoulli measure, we know that the shift map is ergodic with respect to this measure, and the sequence of empirical measures of almost every point is convergent to the Bernoulli measure, hence shift is statistical with respect to this measure. But if the reference measure is taken to be the Dirac mass on the point $\omega$ introduced in this example, then almost every point with respect to this measure would be non-statistical, and hence shift map is non-statistical with respect to $\delta_\omega$. 
\end{remark}

Let us prepare for the proof of Theorem \ref{thm.main.essential}. First, we need to quantify how different are the statistical behaviors of two arbitrary maps $h,g\in \Lambda$ for iterations larger than a fixed number $N\in \mathbb{N}$. To this aim, we propose the following map $\Delta_N$ that associates to a couple of maps $h,g\in\Lambda$ a non-negative real number:

$$ \Delta_{N}(h,g):=\int_{X}{\ \sup_{N\leq n,m}d_w(e_{n}^{h}(x),e_{m}^{g}(x)})\ d\mu,$$
which can be interpreted as the average of the maximum difference between statistical behaviors that the orbit of a point can display under iterations of $h$ and $g$ for iterations larger than $N$. Note that $\Delta_N$ is not a distance. In particular if $f$ is a non-statistical map, then $\Delta_N^e(f,f)$ is uniformly positive for every $N$: 
\begin{lemma}\label{lemma unif} 
A map $f$ is non-statistical if and only if there is a real number $d>0$ such that for each $N\in\mathbb{N}$ we have $\Delta_{N}(f,f)>d$.
\end{lemma}
\begin{proof}
Let $f$ be a non-statistical map, and let $x\in X$ be a non-statistical point. Since the sequence of empirical measures of this point does not converge, 
$$d_x:=\inf_{N>0}\sup_{N\leq n,m}d_w(e_{n}^f(x),e_{m}^f(x))>0$$
By definition, the set of non-statistical points has positive measure and $x\mapsto d_x$ is measurable, thus  
$$\Delta_{N}(f,f)=\int_{X} \ \sup_{N\leq n,m}d_w(e_{n}^f(x),e_{m}^f(x)) d\mu \ge \int_X d_x d\mu>0.$$
To prove the other side let $f$ be a map for which the sequence of empirical measures of almost every point converges. So for almost every $x\in X$ we have 
$$\lim_{N\to\infty} \sup_{N\leq n,m}d_w(e_{n}^f(x),e_{m}^f(x))=0.$$
Since the distance between empirical measures is bounded, we can then use Lebesgue dominated convergence theorem, to conclude 
\begin{align*}
\lim_{N\to \infty} \Delta_{N}(f,f)&=\lim_{N\to \infty}\int_{X}\sup_{N\leq n,m}d_w(e_{n}^f(x),e_{m}^f(x)) d\mu\\
&=\int_{X}\lim_{N\to\infty} \sup_{N\leq n,m}d_w(e_{n}^f(x),e_{m}^f(x))d\mu=0.
\end{align*}
This finishes the proof.
\end{proof}

We recall that $\Lambda$ is a subset of self-mappings of $X$ endowed with a topology finer than $C^0$ topology. Now we want to quantify the difference in the statistical behaviors of maps converging to $f\in \Lambda$. For this purpose, we introduce the following definition:

\begin{definition}
 The \textit{amplitude of statistical divergence} of a map $f\in\Lambda$ is defined as below
\begin{equation*}
    \Delta(f):=\lim_{N\rightarrow \infty} \limsup_{\substack{{h,g\to f}\\ {h,g\in \Lambda}}}\Delta_N(h,g).
\end{equation*}
\end{definition}
Observe that if $\Delta$ is positive at $f$ then the asymptotic behaviors of nearby maps are very sensitive to perturbations of $f$ and so the map $f$ is unstable from the statistical viewpoint. 
\begin{definition}
A map $f\in \Lambda$ is \textit{statistically unstable with respect to $\Lambda$ } if $\Delta(f)>0$.
\end{definition}
\begin{remark}
    Observe that this definition is the same as Definition \ref{def.stat.uns} in the introduction. 
\end{remark}
\begin{example}
Let $X$ be the circle $\mathbb{S}^{1}$, $\mu$ the normalized Lebesgue measure and $\Lambda$ be the set of rotations of the circle. The map $f=Id_{\mathbb{S}^1}$ is statistically unstable since the empirical measures of any point are the Dirac mass at that point, but for any arbitrarily close irrational rotation the sequence of empirical measures converges to the Lebesgue measure.   
\end{example}

\begin{example}
Let $X$ be the Riemann sphere $\hat{\mathbb{C}}$ and $\mu$ its normalized Lebesgue measure. Consider the set of quadratic maps:
\begin{align*}
    \Lambda=\{f_{c}:\hat{\mathbb{C}} \to \hat{\mathbb{C}}| f_{c}(x)=x^2+c \ \ \mathrm{for}\ \  x\in\mathbb{C}, \ \ f_{c}(\infty)=\infty\}.
\end{align*}
    
The map $f_{\frac{1}{4}}$ has a fixed point at $x=\frac{1}{2}$ which attracts the points in a non-empty open set $U$. For any $\epsilon > 0 $, the map $f_{\frac{1}{4}+\epsilon}$ has a different dynamics: almost every point goes to infinity. So for any $\epsilon>0$ and every point $x\in U$ the sequence of empirical measures converges to $\delta_{\frac{1}{2}}$ and $\delta_{\infty}$ under iterating the maps $f_{\frac{1}{4}}$ and $f_{\frac{1}{4}+\epsilon}$ respectively.

Hence the supremum in the definition of $\Delta_N$ is at least $d_w(\delta_{\frac{1}{2}},\delta_{\infty})$ for almost every point. So $\Delta_{N}(f_{\frac{1}{4}},f_{\frac{1}{4}+\epsilon})>\mu(U)d_w(\delta_{\frac{1}{2}},\delta_{\infty})$, which is independent of $\epsilon$ and $N$. According to the definition, this means that $f_{\frac{1}{4}}$ is statistically unstable. 
\end{example}

\begin{proposition}\label{prop Delta^e u.s.c}
The map $\Delta:\Lambda\to \mathbb{R}$ is upper semi-continuous. 
\end{proposition}

\begin{proof}
Let $\{f_k\}_k $ be a sequence of maps converging to $f$. For each $k$ we can find a natural number $N_k$ and two maps $g_k$ and $h_k$ near $f_k$ such that 
$$|\Delta(f_k)-\Delta_{N_k}(g_k,h_k)|<\frac{1}{k}.$$
Note that we can choose the sequence $\{N_k\}_k$ such that it converges to infinity and also the sequence of maps $\{g_k\}_k$ and $\{h_k\}_k$ such that both converge to $f$. So we obtain 

\begin{equation*}
    \Delta(f)\geq \limsup_{k}\Delta(f_k),
\end{equation*}
and this implies the upper semi-continuity of $\Delta$.
\end{proof}

\section{Proof of Theorem \ref{thm.main.essential}}
Now we are ready to prove Theorem \ref{thm.main.essential}:

\begin{proof}[Proof of Theorem \ref{thm.main.essential}]
Since by Proposition \ref{prop Delta^e u.s.c} the map $\Delta$ is an upper semi-continuous map, there is a Baire generic set $\mathcal{G}$ on which the map $\Delta$ is continuous. For a map $f \in \mathcal{G}$, which is also in the interior of statistically unstable maps, there exists a neighborhood $\mathcal{U}_f$ around $f$ on which $\Delta$ is uniformly positive:
\begin{equation*}
    \exists d>0 \quad  s.t. \quad  \forall g \in \mathcal{U}_f,\quad  \Delta(g)>d. 
\end{equation*}

Now we construct a sequence of open and dense subsets in ${\mathcal U}_f$ such that any map in the intersection of these sets is non-statistical. And then we can conclude that non-statistical maps are Baire generic in ${\mathcal U}_f$.

To construct such open and dense sets we need to show a semi-continuity property of the map $\Delta_N$:
\begin{lemma}
The map $\Delta_{N}$ is lower semi-continuous.
\end{lemma}
\begin{proof}
We recall that for $h,g\in \Lambda$ 
$$ \Delta_{N}(h,g):=\int_{X}{\ \sup_{N\leq n,m}d_w(e_{n}^{h}(x),e_{m}^{g}(x)})\ d\mu,$$

Now note that 
$$ \Delta_{N}(h,g)= \sup_{N\leq M}\int_{X} {\ \sup_{N\leq n,m \leq M}d_w(e_{n}^{h}(x),e_{m}^{g}(x))\ d\mu}, $$
and we know that $\int_{X} {\ \sup_{N\leq n,m \leq M}d_w(e_{n}^{h}(x),e_{m}^{g}(x))\ d\mu}$ is continuous with respect to $h$ and $g$ and supremum of a sequence of continuous functions is lower semi-continuous so we are done. 
\end{proof}
Next, we show that for any $N\in\mathbb N$, the set 
$$E(N):=\{h\in {\mathcal U}_f|\Delta_N(h,h)>\frac{d}{3}\},$$
is an open and dense subset of ${\mathcal U}_f$ and moreover every map in the intersection $\bigcap_N E(N)$ is non-statistical.

The openness of $E(N)$ is guaranteed by lower semi-continuity of $\Delta_N.$ Now we prove the denseness of $E(N)$. For any arbitrary map $h\in {\mathcal U}_f\subset \Lambda$, take a neighborhood $V_h$ such that for any map $g\in V_h$ and any $x\in X$:
\begin{equation}\label{p.w. in.eq d/3}
    d_w(e^h_{N}(x),e^g_{N}(x))<\frac{d}{3}.
\end{equation}
This is possible because $N$ is fixed, $e^g_{N}$ depends continuously on $g$ and $X$ is compact. Now since $\Delta(h)>d$, we can choose $g_1,g_2\in V_h$ such that for some integer $M>N$ it holds true that 
\begin{equation}
    \Delta_M(g_1,g_2)>d.
\end{equation}
But since $\Delta_N(g_1,g_2)$ is decreasing in $N$ we obtain:

\begin{equation}\label{p.w. in.eq d}
    \Delta_N(g_1,g_2)>d.
\end{equation}
Now note that for each $x\in X$ we have 
\begin{multline*}
    \quad\quad\quad\quad\quad\quad\quad \sup_{N\leq n,m} d_w(e_{n}^{g_1}(x),e_{m}^{g_2}(x))\leq \\
\sup_{N\leq n,m}d_w(e_{n}^{g_1}(x),e_{m}^{g_1}(x)) + \sup_{N\leq n,m}d_w(e_{n}^{g_2}(x),e_{m}^{g_2}(x))+ d_w(e_{N}^{g_1}(x),e_{N}^{g_2}(x)), 
\end{multline*}

and hence after integrating with respect to $\mu$ and using inequality \ref{p.w. in.eq d/3} we obtain:

$$\Delta_N(g_1,g_2)\leq \Delta_N(g_1,g_1)+\Delta_N(g_2,g_2)+ \frac{d}{3}.$$
Now using inequality \ref{p.w. in.eq d} we conclude that at least one of the maps $g_1$ and $g_2$ is inside $E(N)$, and then recalling that $h$ was chosen arbitrarily in ${\mathcal U}_f$ and $V_h$ arbitrary small, we conclude that $E(N)$ is dense in ${\mathcal U}_f$.
  
Observe that Lemma \ref{lemma unif} implies that any map $h$ in the set $\bigcap_{N=1}^\infty E(N)$, which is a Baire generic set inside ${\mathcal U}_f$, is non-statistical. So the non-statistical maps are generic in some open neighborhoods of $f$ and $f$ is an arbitrary map in the interior of statistically unstable maps intersected by the generic set $\mathcal G$. This implies that non-statistical maps are indeed a generic subset of the interior of statistically unstable maps as well. 

\end{proof}

\section{Proof of Theorem \ref{theorem.generic.phase}}\label{Sec.gen.non-stat.phase}
In this section, we show how our results for the Baire genericity of non-statistical maps in a subset of dynamics, can be used to obtain the Baire genericity of non-statistical points in the phase space for a single map. \\
Suppose $f:X\rightarrow X$ is a continuous map of a compact metric space $X$. Using this map, we will construct a family of dynamical systems parameterized by points of the space $X$. Add an auxiliary point $\hat{x}$ to the space $X$ to define a new compact metric space:
$$\hat{X}:=\{\hat{x}\}\cup X.$$
We define the distance between $\hat{x}$ and any point in $X$ to be equal to $1$. For each point $x\in X$ define the map $F_x:\hat{X}\rightarrow\hat{X}$ such that $F_x(\hat{x})=x$ and $F_x|_X=f$. These maps will form a family $\Lambda:=\{F_x\}_{x\in X}$ of maps on $\hat{X}$, parameterized by points in $X$. Take the reference probability measure $\mu$ on the space $\hat{X}$ to be Dirac mass $\delta_{\hat{x}}$ on the auxiliary point $\hat{x}$. The push forward of this measure under any map $F_x$ is the Dirac mass on the point $x$. After the first iteration, it is indeed the map $f$ that is iterated, and so this reference probability measure captures the statistical behavior of the $f$-orbit of $x$ for the map $F_x$. We induce the topology of the space $X$ to the set $\Lambda$. We can apply our results in this case. 
Observe that a map $F_x\in \Lambda$ is non-statistical if and only if the point $x\in X$ is non-statistical.\\
\begin{definition}
We say a point $x\in X$ is statistically unstable if the map $F_x$ is statistically unstable as an element of $\{F_x\}_{x\in X}$. 
\end{definition}

Now we are ready to prove Theorem \ref{theorem.generic.phase}:

\begin{proof}[Proof of Theorem \ref{theorem.generic.phase}]
By definition, every map in $\Lambda=\{F_x\}_{x\in X}$ is statistically unstable. Using Theorem \ref{thm.main.essential} we conclude that a Baire generic map $F_x$ is non-statistical. This indeed implies that the point $x$ is non-statistical for the map $f$. \end{proof}

In \cite{Carvalho2021SensitivityAH} sufficient conditions for an observable to have non-convergent Ces\'aro averages on a Baire generic set in the phase space are introduced. Theorem \ref{theorem.generic.phase} is similar to their result. They also define a notion called sensitivity that is similar to the notion of statistical instability. It can be checked that their corollaries can be obtained from Theorem \ref{theorem.generic.phase} as well. We should mention that their result covers the case of an arbitrary sequence of bounded observables, and not only observables coming from Ces\'aro averages. Let us now prove corollary \ref{cor.minimal.not.erg}:
\begin{proof}[Proof of Corollary \ref{cor.minimal.not.erg}]
   As the map $f$ is not uniquely ergodic, it has at least two different ergodic measures $\nu_1$ and $\nu_2$. As the map is minimal, the support of any invariant measure is the whole space. So for any $x\in X$, and any neighborhood $U$ of $x$, there is two point $y_1$ and $y_2$ in $U$ such that 
   $$\lim_{n\rightarrow \infty} e_n^f(y_i)=\nu_i, \quad i\in \{1,2\}.$$ This will imply that for the family of maps $\Lambda=\{F_x|x\in X\}$ defined as above, we have 
   $$\Delta_\Lambda(F_x)\geq d_w(\nu_1,\nu_2),$$
   and hence $F_x$ is statistically unstable. This will imply that $x$ is statistically unstable, so we can use Theorem \ref{theorem.generic.phase}.
\end{proof}

\section{Proof of Theorem \ref{Anosov-Katok example}}\label{sec.Anosov-Katok}

Let us first give an overview of the proof. In Lemma \ref{lemma-AK} we use the Anosov-Katok method to show that the maps in $\mathcal{AK}$ are statistically unstable. Indeed we show that near any map and for any boundary component, we can find another map for which most of the points spend most of the time close to that boundary component. This will give us statistical instability and up to here, we can use Theorem \ref{thm.main.essential} to obtain Baire genericity of non-statistical maps in $\mathcal{AK}$. But we go further to use the theorem of Fayad and Katok \cite[Theorem 3.3]{FK04} which describes Baire generic behavior of conservative Anosov-Katok maps of the annulus to deduce Baire generic behavior in the dissipative setting and obtain maximally oscillation. Here we need to use some new definitions and results developed in the version of convergence in law, which can be found in Appendix \ref{bif-to}. Let us recall some of the definitions.\\

For each positive integer $n$, $e_n^f$ defines a map from the phase space to its set of probability measures:
$$e_n^f:X\rightarrow \mathcal M_1(X).$$
We can investigate how the empirical measures $e^f_n(x)$ are distributed in $\mathcal M_1(X)$ with respect to the reference measure $\mu$ on $X$ and what is the asymptotic behavior of these distributions. To this aim, we can push forward the measure $\mu$ to the set of probability measures on $X$ using the map $e^f_n$: 
$$\hat{e}_n(f):=(e_n^f)_*(\mu).$$ 
The measure $\hat{e}_n(f)$ is a probability measure on the space of probability measures on $X$. We denote the space of probability measures on the space of probability measures by $\mathcal M_1(\mathcal M _1(X))$. We denote the Wasserstein metric on this space by $\hat{d}$. Note that the compactness of $X$ implies the compactness of $\mathcal M_1(X)$ and hence the compactness of $\mathcal M_1(\mathcal M_1(X))$. So the sequence $\{\hat{e}_n(f)\}_{n\in \mathbb{N}}$ lives in a compact space and has one or possibly more than one accumulation points.

\begin{definition}
 For a map $f\in\Lambda$ and a probability measure $\hat{\nu}\in \mathcal M_1(\mathcal M_1(X))$, we say $f$ \emph{statistically bifurcates toward $\hat{\nu}$ through perturbations in $\Lambda$}, if there is a sequence of maps $\{f_k\}_k $ in $\Lambda$ converging to $f$ and a sequence of natural numbers $\{n_k\}_k$ converging to infinity such that the sequence $\{\hat{e}_{n_k}(f_k)\}_k$ converges to $\hat{\nu}\in\mathcal M_1(\mathcal M_1(X))$. 
\end{definition}

  We remind that for any measure $\nu\in\mathcal M_1(X)$, the Dirac mass on $\nu$, which is an element of $\mathcal M_1(\mathcal M_1(X))$, is denoted by $\hat{\delta}_\nu$. \\

Now we are ready to begin proving Theorem \ref{Anosov-Katok example}. We start with the following lemma:
\begin{lemma}\label{lemma-AK}
Let $C$ be one of the connected components of the boundary of $\mathbb{A}$, and $f$ be an arbitrary map in $\mathcal{AK}$. Then there is a measure $\nu$ supported on the set $C$ such that $f$ statistically bifurcates towards $\hat{\delta}_{\nu}$.
\end{lemma}

\begin{proof}
Since $f$ can be approximated by maps that are conjugate to rotation, there is a rational number $\frac{p}{q}$ and a $C^r$ diffeomorphism $h$ such that the map $h^{-1}R_{\frac{p}{q}}h$ is close to $f$.  
Let  
$$B:=[r_1,r_2]\times \mathbb R/ \mathbb Z,$$
for some distinct real numbers $r_1,r_2\in (0,1)$. Take a real number $\theta>0$ and define 
$$B_1:= [r_1,r_2]\times [0,\theta)$$
and
$$B_2:= [r_1,r_2]\times [\theta,1). $$

\begin{sublemma}\label{sublemma}
For any $\sigma< 1$ close to $1$ and $\epsilon>0$ small, there is a map $\hat{g} \in\mathrm{Diff}_{+}^r(\mathbb{A})$ with the following properties:
\begin{itemize}
    \item $\hat{g}$ is identity on a neighborhood of the set $C$,
    \item $Leb(\hat{g}(B_{1}))>\sigma$,      
    \item $\hat g(B_2)$ is included in the $\epsilon$-neighborhood $N_\epsilon(C)$ of $C$.
\end{itemize}
\end{sublemma}
\begin{figure}
    \centering
    \includegraphics[scale=0.15]{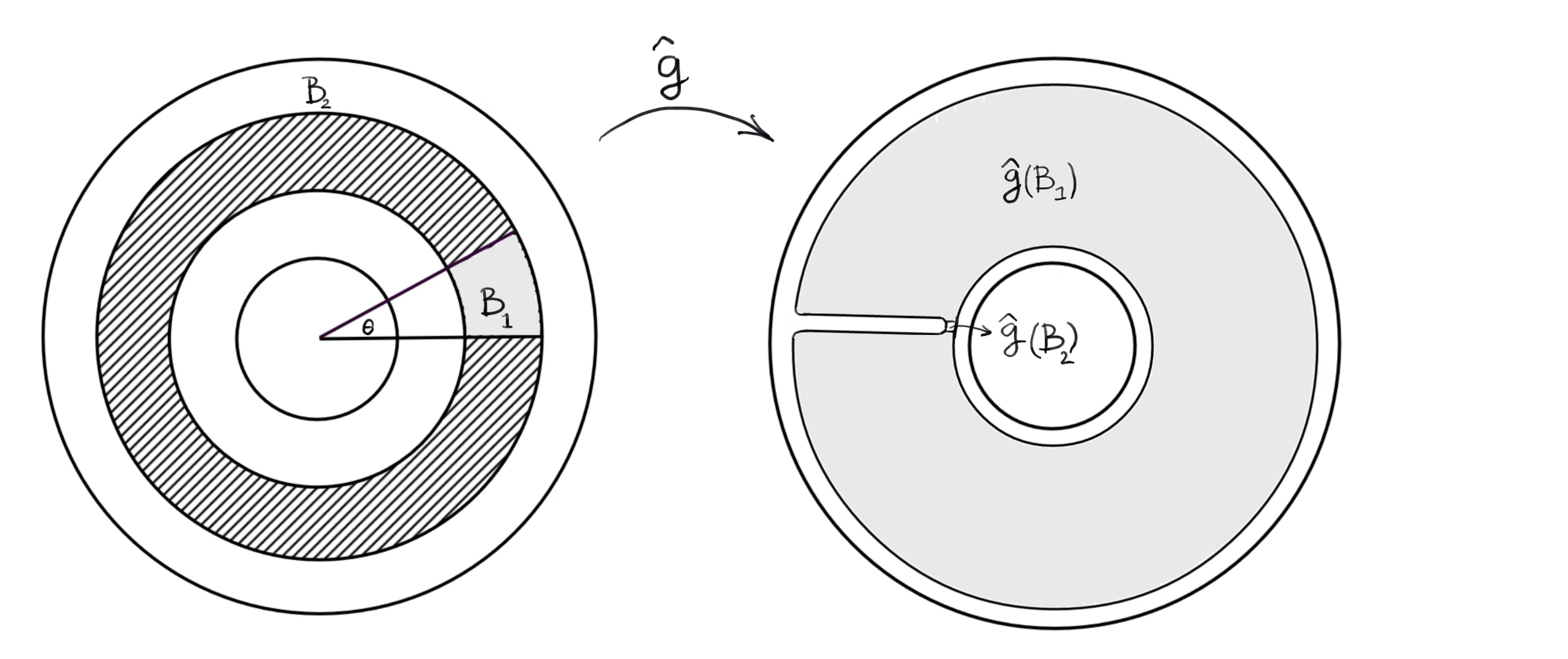}
    \caption{The map $\hat{g}$}
    \label{Anosov-Katok-pic}
\end{figure}
\begin{proof}
Using bump functions, we construct a map $\hat g$ as depicted in Figure \ref{Anosov-Katok-pic}. The technical details are left to the reader.
\end{proof}

Now let $\hat{g}$ be a map found in the sublemma. This map can be lifted using the covering map $\pi:\mathbb{A}\to \mathbb{A}$, $\pi(r,\theta)=(r,q\theta)$. Let $g$ be the lift of $\hat{g}$ which is identity around the set $C$. The diffeomorphism $g$ has similar properties:
\begin{itemize}
    \item ${g}$ is identity on a neighborhood of the set $C$,
    \item $Leb({g}(\pi^{-1}(B_{1})))>\sigma$,
    \item ${g}(\pi^{-1}(B_{2}))\subset N_\epsilon(C)$.
\end{itemize}
Now note that $g$ commutes with $R_{\frac{p}{q}}$ so:
\begin{equation*}
    h\circ g \circ R_{\frac{p}{q}} \circ g^{-1} \circ h^{-1}=h\circ  R_{\frac{p}{q}}  \circ h^{-1}
\end{equation*}
Choose $\alpha'$ irrational and small enough so that $h\circ g \circ R_{\alpha'} \circ g^{-1} \circ h^{-1}$ is arbitrary close to $h\circ  R_{\frac{p}{q}}  \circ h^{-1} $. Indeed $h\circ g$ is $C^r$ and the map sending $\alpha$ to  $h\circ g \circ R_{\alpha} \circ g^{-1} \circ h^{-1}$ is hence continuous. Since $\alpha'$ is irrational, the orbit closure of any point in $\mathbb{A}$ under iterating the map $h\circ g \circ R_{\alpha'} \circ g^{-1} \circ h^{-1}=:f'$ is $h\circ g$-image of the orbit closure of a point under iterating the map $R_{\alpha'}$, which is a circle in $\mathbb{A} $. So for any point $x\in h\circ g(\pi^{-1}(B))$, the orbit closure of $x$ is the $h\circ g$ image of a vertical circle $C'$ contained in $\pi^{-1}(B)$. The map $h\circ g \circ R_{\alpha'} \circ g^{-1} \circ h^{-1}|_{h\circ g(C')}$ is conjugate to $R_{\alpha'}|_{C'}$. 
Now note that if $R_{\alpha'}^{n}\circ g^{-1}\circ h^{-1}(x) \in \pi^{-1}(B_{i})$ then $(f')^{n}(x)\in h\circ g(\pi^{-1}(B_{i}))$ for $i\in \{1,2\}$. The orbit of each point in $C'$, in average, spends ${\theta}$ portion of times in $\pi^{-1}(B_1)$ and ${1-\theta}$ portion of times in $\pi^{-1}(B_2)$. So the orbit of the point $x$ spends ${\theta}$ portion of times in $h\circ g(\pi^{-1}(B_1))$ and ${1-\theta}$ portion of times in $h\circ g(\pi^{-1}(B_2))$. By choosing $\theta$ and $\epsilon$ sufficiently small, we can guarantee that the asymptotic averages of any point in $h\circ g(\pi^{-1}(B_1))$ is arbitrarily close to a measure $\nu_{f'}$ supported on $C$ which is the pushforward of the Lebesgue measure on the boundary component $C$ by the map $h\circ g$. This is because the image of $\pi^{-1}(B_1)$ under the map $g$ can be so close to $C$ that after iteration by h, it still remains arbitrarily close to $C$. Since $h$ is smooth and sends positive measure sets to positive measure sets, it can be shown that the sets with large measure are sent to large measure sets and hence, If $\sigma$ is chosen sufficiently close to one, then map $g$ is so that  $Leb(h\circ g(\pi^{-1}(B_1)))$ is sufficiently close to one and so since for a large number $n$ the $n^{th}$ empirical measures of points in $h\circ g(\pi^{-1}(B_1))$ is close to $\nu_{h\circ g}$ and so the measure $\hat{e}^{f'}_n$ is close to $\hat{\delta}_{\nu_{f'}}$. Now taking $\nu$ as any accumulation point of measures like $\nu_{f'}$ where $f'$ approaches $f$, according to the definition we can see that $f$ statistically bifurcates toward $\hat{\delta}_{\nu}$ and so we are done.
\end{proof}
\begin{proof}[Proof of Theorem \ref{Anosov-Katok example}]
 Lemma \ref{lemma-AK} implies that for any map $f\in \mathcal{AK}^r$ there are two measures $\nu_{1,f}$ and $\nu_{2,f}$ which are supported on different connected components of the boundary of $\mathbb{A}$ such that $f$ statistically bifurcates toward both $\hat{\delta}_{\nu_{1,f}}$ and $\hat{\delta}_{\nu_{2,f}}$. Now using Theorem \ref{generic stat behav} we conclude that for a generic map $f\in\mathcal{AK}^r$ and for almost every point $x$ in the phase space $X$, the set of accumulation points of the sequence $\{{e}_n^f(x)\}_n$ contains at least two measures ${\nu_{1,f}}$ and ${\nu_{2,f}}$. We are going to show that generically these two measures are the only ergodic invariant measures of the map $f$ and the empirical measures of almost every point accumulate to any convex combination of these two ergodic measures (which is indeed the whole space of invariant measures) and hence $f$ is maximally oscillating. Now approximate $f$ by a map like $h\circ g\circ R_{\alpha'}\circ g^{-1}\circ h^{-1}$ coming from Lemma \ref{lemma-AK}. Using the theorem of Fayad and Katok, we know that the map $R_{\alpha'}$ can be approximated in $C^\infty$-topology (and hence in $C^r$-topology) by a map $e\in \mathcal{AK}_{vol}^\infty$ which has only three ergodic measures, two one dimensional Lebesgue measures on the boundaries and the volume measure of the annulus. The map $h\circ g\circ e \circ g^{-1} \circ h^{-1}$ is close to the initial map $f$ and has only three ergodic invariant measures which are the push forward of three ergodic measures of $e$ by the map $h\circ g$. Note that if in Sublemma \ref{sublemma} the numbers $r_1$ and $r_2$ are chosen close to $0$ and $1$ then the set $B_2$ has measure close to one. In this case, observe that the pushforward of the volume measure by the map $h\circ g$ is a measure that is close to the pushforward of the one-dimensional Lebesgue measure of one of the boundary components (which is denoted by $C$ in the lemma). Hence the set of invariant measures for the map $h\circ g\circ e\circ g^{-1}\circ h^{-1}$ is a triangle that has two of its vertices very close to each other. We know that the map sending the dynamics to its set of invariant measures is upper semi-continuous (see \ref{semi-cont inv meas}) and hence it is continuous for maps in a Baire generic set. Hence we can assume that $f$ is a continuity point of this mapping and the set of invariant measures of $f$ is in an arbitrary neighborhood of a triangle which is arbitrarily close to a segment. So the set of invariant measures of $f$ is a segment. But this means that $f$ has exactly two ergodic invariant measures. Note that these measures are supported by different boundary components of the annulus. So on each boundary component of the annulus, there is only one ergodic measure, and hence any invariant measures on one of the boundary components are equal to the corresponding ergodic measure on that component. So two measures $\nu_{1,f}$ and $\nu_{2,f}$ toward which $f$ statistically bifurcates are exactly two ergodic measures of $f$. Moreover, since these two measures are in the accumulation points of the sequence of empirical measures for almost every point, and the set of invariant measures is the line segment between these two measures, the sequence of empirical measures of almost every point has to accumulate to every point in this line segment and this finishes the proof. 
\end{proof}

\appendix
\section{The $L^1$ convergence version}
We recall the $n^{\text{th}}$ empirical function of a map $f$ is the map $e^f_n:X\to \mathcal M_1(X)$ sending a point $x\in X$ to the $n^{\text{th}}$ empirical measure $e^f_n(x)$. We are going to study the $L^1$ (non-) convergence of the sequence of empirical functions. For this purpose, we need to give some definitions:

Let us denote the space of Borel measurable maps from $X$ to $\mathcal M_1(X) $ by $L^1(X,\mathcal M_1(X))$. Note that since the empirical functions are continuous maps with respect to $x$, they are elements of $L^1(X,\mathcal M_1(X))$.
We define a metric on this space where the distance between two elements $e,e'\in L^1(X,\mathcal M_1(X))$ is defined as follows:
\begin{equation}\label{L^1 distance}
d_{L^1}(e,e')=\int _X d_w(e(x),e'(x))d\mu.
\end{equation} 
Let us study the convergence of the sequence of empirical functions with respect to this metric:
\begin{definition}
We say a map $f$ is \textit{$L^1$ non-statistical} if the sequence of maps $e^f_{n}:X\to \mathcal M_1(X)$ is not convergent for the $L^1$ topology:

\begin{equation*}
    \limsup_{m,n\to\infty}d_{L^1}(e^f_{n},e^f_{m})>0.
\end{equation*}

\end{definition}

In the following we introduce a condition, that if it is satisfied by a Baire space of dynamics $\Lambda$, then we can conclude the existence of $L^1$ non-statistical maps within $\Lambda$. 
Let us first, quantify the extent to which the statistical behavior of a map $f\in \Lambda$ can be changed by small perturbations. To this aim, we propose the following definition:
\begin{definition}
 The \emph{amplitude of $L^1$ statistical divergence w.r.t. $\Lambda$} of a map $f\in\Lambda$ is a real-valued non-negative mapping which is defined as follows 
\begin{equation*}
    \Delta^1(f):= \limsup_{\substack{{h,g\to f,m,n\to \infty}\\ {h,g\in \Lambda}}} d_{L^1}(e^h_{n},e^g_{m}).
\end{equation*}
\end{definition}
Note that the definition of $\Delta^1$ depends also on the set $\Lambda$, and not only on the map $f$. However, for the sake of simplicity, we hid this in the notation. 

Observe that if $\Delta^1$ is positive at $f$ then the asymptotic behaviors of dynamics close to $f$ are very sensitive to perturbations and in this sense, the map $f$ is statistically unstable. We introduce the following definition: 
\begin{definition}
A map $f\in \Lambda$ is $L^1$ \textit{statistically unstable with respect to $\Lambda$ } if its amplitude of $L^1$ statistical divergence w.r.t. $\Lambda$ is positive:
$$\Delta^1(f)>0.$$
\end{definition}
\begin{example}\label{exp-IdS1}
Suppose $\Lambda$ is the set of rigid rotations on $\mathbb{S}^1$. The identity map $Id_{\mathbb{S}^1}\in\Lambda$ is $L^1$ statistically unstable, since the empirical measures of all of its points are atomic whereas we can approach the map $Id_{\mathbb{S}^1}$ by irrational rotations, and the empirical measures of any point are close to the Lebesgue measure for large enough iterations. 
\end{example}
Now we want to investigate the relationship between $L^1$ statistical instability and the existence of $L^1$ non-statistical maps. It is clear that if a map $f$ is $L^1$ non-statistical then $\Delta^1(f) >0$ and so $f$ is $L^1$ statistically unstable, but the existence of a $L^1 $ statistically unstable map does not necessarily imply the existence of $L^1$ non-statistical maps (see Example \ref{exp-IdS1} ). However, if the interior of $L^1$ statistically unstable maps is non-empty, then the existence of plenty of $L^1$ non-statistical maps is guaranteed. 

\begin{customthm}{D}\label{thm.main.L^1}
The $L^1$ non-statistical maps form a Baire generic subset of the interior of $L^1$ statistically unstable maps. 
\end{customthm}
\begin{remark}
    Theorem \ref{thm.main.L^1} is the counterpart of Theorem \ref{thm.main.essential} in the $L^1$ convergence version. 
\end{remark}

The proofs of Theorem \ref{thm.main.L^1} uses the following lemma:
\begin{lemma}\label{lem.l.s.c.Delta'}
The map $\Delta^1$ is upper semi-continuous.
\end{lemma}

\begin{proof}

Let $\{f_k\}_k $ be a sequence of maps converging to $f$. For each $k$ we can find two natural numbers $n_k$ and $m_k$ and two maps $g_k$ and $h_k$ near $f_k$ such that 
$$|\Delta^1(f_k)-d_{L^1}(e_{n_k}^{g_k},e_{m_k}^{h_k})|<\frac{1}{k}.$$
Note that we can choose the sequences $\{n_k\}_k$ and $\{m_k\}_k$ both converging to infinity and also the sequence of maps $\{g_k\}_k$ and $\{h_k\}_k$ both converging to $f$. So we obtain 

\begin{equation*}
    \Delta^1(f)\geq \limsup_{k}\Delta^1(f_k),
\end{equation*}
and this implies the upper semi-continuity of $\Delta^1$.
\end{proof}

\begin{proof}[Proof of Theorem \ref{thm.main.L^1}]
Since $\Delta^1$ is upper semi-continuous, there is a generic subset $\mathcal{G}\subset \Lambda$ on which the map $\Delta^1$ is continuous. For a map $f \in \mathcal{G}$ which is also in the interior of $L^1$ statistically unstable maps, there exists a neighborhood $\mathcal{U}_f\subset \Lambda$ around $f$ on which $\Delta^1$ is uniformly positive:
\begin{equation}\label{equ unif}
    \exists d>0 \quad  s.t. \quad  \forall g \in \mathcal U_f,\quad  \Delta^1(g)>d. 
\end{equation}

Now we construct a sequence of open and dense subsets in ${\mathcal U}_f$ such that any map in the intersections of these sets is $L^1$ non-statistical. This will imply that $L^1$ non-statistical maps are Baire generic in ${\mathcal U}_f$ and hence the $L^1$ non-statistical maps are locally generic in the interior of $L^1$ statistically unstable maps. 

To construct such open and dense sets, first note that the function $\Delta^1$ can be written as 
\begin{equation*}
\Delta^1(f)=\limsup_{g,h \to f,N \to \infty} \Delta ^{1}_{N}(h,g),
    \end{equation*}
    where $\Delta^1_{N}(h,g)=\sup_{i,j\ge N}d_{L^1}(e^h_{i},e^g_{j}).$

\begin{claim}
The map $\Delta^1_{N}$ is lower semi-continuous.
\end{claim}

\begin{proof}
Note that 
$$\Delta^1_{N}(h,g)=\sup_{M\ge N}\{\sup_{N \le i,j \le M}d_{L^1}(e^h_{i},e^g_{j})\}.$$

But $\sup_{N \le i,j \le M}d_{L^1}(e^h_{i},e^g_{j})$ is continuous with respect to $(h,g)$, and supremum of a sequence of continuous functions is lower semi-continuous. 
\end{proof}

Next, we show that for any $N\in\mathbb N$, the set 
$$E(N):=\{h\in {\mathcal U}_f|\Delta^1_N(h,h)>\frac{d}{3}\},$$
is an open and dense subset of ${\mathcal U}_f$ and moreover every map in the intersection $\bigcap_N E(N)$ is $L^1$ non-statistical.

The openness of $E(N)$ is guaranteed by lower semi-continuity of $\Delta^1_N.$ Now we prove the denseness of $E(N)$. For any arbitrary map $h\in {\mathcal U}_f$, take a neighborhood $V_h$ such that for any map $g\in V_h$ we have:
\begin{equation}\label{in.eq d/3}
    d_{L^1}(e^h_{N},e^g_{N})<\frac{d}{3}.
\end{equation}
This is possible because $N$ is fixed and $e^g_{N}$ depends continuously on $g$. By \ref{equ unif} we know that $\Delta^1(h)>d$, and so we can choose $g_1,g_2 \in V_h$ such that for some integers $n,m>N$ it holds true that 
\begin{equation}\label{in.eq d}
    d_{L^1}(e^{g_1}_{n},e^{g_2}_{m})>d.
\end{equation}
But note that 
\begin{equation*}
    d_{L^1}(e^{g_1}_{n},e^{g_2}_{m})\le d_{L^1}(e^{g_1}_{n},e^{g_1}_{N})+d_{L^1}(e^{g_1}_{N},e^{g_2}_{N})+d_{L^1}(e^{g_2}_{N},e^{g_2}_{m}).
\end{equation*}
  Inequalities \ref{in.eq d/3} and \ref{in.eq d} imply that either $$d_{L^1}(e^{g_1}_{n},e^{g_1}_{N})>\frac{d}{3} \text{ or } d_{L^1}(e^{g_2}_{m},e^{g_2}_{N})>\frac{d}{3}.$$
  So at least one of the maps ${g_1}$ or ${g_2}$ is inside $E(N)$, and then recalling that $h$ was chosen arbitrarily in ${\mathcal U}_f$ and $V_h$, we conclude that $E(N)$ is dense in ${\mathcal U}_f$.
    
  Now observe that for any map $h\in\bigcap_{N=1}^\infty E(N)$, the sequence of empirical functions is not a Cauchy sequence and hence $h$ is $L^1$ non-statistical. The set $\bigcap_{N=1}^\infty E(N)$ is a Baire generic set in the open neighborhood $\mathcal U_f$, so the $L^1$ non-statistical maps are generic in the set $\mathcal U_f$. Considering the fact that $f$ is an arbitrary map in the generic set $\mathcal G$ we can then conclude that $L^1 $ non-statistical maps are indeed a generic subset of the interior of $L^1$ statistically unstable maps. 
\end{proof}

\section{The version of convergence in law}\label{bif-to}
For a dynamical system $f:X\to X$ the map $e^f_n:X\to \mathcal M _1(X)$ associates to each point $x\in X$, its $n^{th}$ empirical measure. Different points usually have different empirical measures. We can investigate how the empirical measures $e^f_n(x)$ are distributed in $\mathcal M_1(X)$ with respect to the reference measure $\mu$ on $X$ and what is the asymptotic behavior of these distributions. To this aim, we can push forward the measure $\mu$ to the set of probability measures on $X$ using the map $e^f_n$: 
$$\hat{e}_n(f):=(e_n^f)_*(\mu).$$ 
The measure $\hat{e}_n(f)$ is a probability measure on the space of probability measures on $X$. We denote the space of probability measures on the space of probability measures by $\mathcal M_1(\mathcal M _1(X))$. We denote the Wasserstein metric on this space by $\hat{d}$. Note that the compactness of $X$ implies the compactness of $\mathcal M_1(X)$ and hence the compactness of $\mathcal M_1(\mathcal M_1(X))$. So the sequence $\{\hat{e}_n(f)\}_{n\in \mathbb{N}}$ lives in a compact space and has one or possibly more than one accumulation points.
\begin{definition}
    A map $f$ is called \textit{statistical in law} if the sequence $\{\hat{e}_n(f)\}_{n\in \mathbb{N}}$ is convergent. 
\end{definition}
\begin{example}
For any $\mu$ preserving map $f:X\to X$, the sequence $\{\hat{e}_n(f)\}_{n\in \mathbb{N}}$ converges to a measure $\hat{\mu}$ which is the ergodic decomposition of the measure $\mu$.
\end{example}
\begin{example}
If $\nu$ is a physical measure for the map $f:X\to X$ whose basin covers $\mu$-almost every point, the sequence $\{\hat{e}_n(f)\}_{n\in \mathbb{N}}$ converges to the Dirac mass concentrated on the point $\mu\in \mathcal M_1(X)$, which we denote by $\hat{\delta}_\mu$.
\end{example}

The following lemma provides some information about the sequence $\{\hat{e}_n(f)\}_{n\in \mathbb{N}}$:
\begin{lemma}\label{lemma.dec.dist}
For any $f\in\Lambda$ and any $n\in\mathbb{N}$ it holds true that 
$$\hat{d}_w(\hat{e}_n^f,\hat{e}_{n+1}^f)< \frac{diam(X)}{n+1},$$
where $diam(X)$ is the diameter of the space $X$. 
\end{lemma}
\begin{proof}
We recall that 
$$\hat{e}_n^f=(e_n^f)_{*}(\mu).$$
First, let us show for any $x\in X$ the following inequality holds true independent of the choice of $f\in\Lambda$: 
$$d_w(e^f_n(x),e^f_{n+1}(x))<\frac{diam(X)}{n+1}.$$
So according to the definition, we should show that
$$\inf_{\gamma\in\pi(e^f_n(x),e^f_{n+1}(x))}\int_{X\times X}d(x,y)d\gamma(x,y)<\frac{diam(X)}{n+1}.$$
Consider the following element of $\pi(e^f_n(x),e^f_{n+1}(x))$:
$$\gamma=\frac{1}{n+1}\Sigma_{0\leq i \leq n-1}\delta_{(f^i(x),f^i(x))}+\frac{1}{n(n+1)}\Sigma_{0\leq i \leq n-1}\delta_{(f^i(x),f^n(x))}.$$
Note that
$$(\pi_1)_*(\gamma)=e_n^f(x)=\frac{1}{n}\Sigma_{0\leq i\leq n-1}\delta_{f^i(x)},$$

$$(\pi_2)_*(\gamma)=e_{n+1}^f(x)=\frac{1}{n+1}\Sigma_{0\leq i\leq n}\delta_{f^i(x)},$$
where $\pi_1$ and $\pi_2$ are the projection on the first and second coordinates. So we have $\gamma\in\pi(e^f_n(x),e^f_{n+1}(x))$ and hence 
\begin{align*}
    {d}_w(e^f_n(x),e^f_{n+1}(x))\leq & \int_{X\times X}d(x,y)d\gamma(x,y)\\
    &=\Sigma_{0\leq i\leq n-1}\frac{1}{n(n+1)}d(f^i(x),f^n(x))\\
    &\leq\frac{diam(X)}{n+1}.
\end{align*}
Now consider the following measure on $\mathcal M_1(\mathcal M_1(X))\times \mathcal M_1(\mathcal M_1(X))$:
$$\hat{\gamma}=\int_X \delta_{(e^f_n(x),e^f_{n+1}(x)}d\mu.$$
Obviously $\hat{\gamma}\in \pi(\hat{e}^f_n,\hat{e}^f_{n+1})$, and so
\begin{align*}
    \hat{d}_w(\hat{e}^f_n,\hat{e}^f_{n+1})
    &\leq \int_{\mathcal M_1(\mathcal M_1(X))\times \mathcal M_1(\mathcal M_1(X))}d(x,y)d\hat{\gamma}(x,y)\\ 
    &=\int_X {d}_w(e^f_n(x),e^f_{n+1}(x))d\mu \le \frac{diam(X)}{n+1}.
\end{align*}  

\end{proof}
Now let $\Lambda$ be a Baire space of self-mappings of $X$ endowed with a topology finer than $C^0$-topology. For each $f\in\Lambda$ the accumulation points of the sequence $\{\hat{e}_n(f)\}_{n\in \mathbb{N}}$ form a compact subset of $\mathcal M_1(\mathcal M_1(X))$ which we denote it by $acc( \{\hat{e}_n(f)\}_{n\in \mathbb{N}})$. This set can vary dramatically by small perturbations of $f$ in $\Lambda$:

\begin{definition}
 For a map $f\in\Lambda$ and a probability measure $\hat{\nu}\in \mathcal M_1(\mathcal M_1(X))$, we say $f$ \emph{statistically bifurcates toward $\hat{\nu}$ through perturbations in $\Lambda$}, if there is a sequence of maps $\{f_k\}_k $ in $\Lambda$ converging to $f$ and a sequence of natural numbers $\{n_k\}_k$ converging to infinity such that the sequence $\{\hat{e}_{n_k}(f_k)\}_k$ converges to $\hat{\nu}\in\mathcal M_1(\mathcal M_1(X))$. 
\end{definition}

For the sake of simplicity, when the space $\Lambda$ in which we are allowed to perturb our dynamics is fixed, we say $f$ statistically bifurcates toward $\hat{\nu}$.

For any map $f\in\Lambda$, by $\mathcal{B}_{\Lambda,f}$ we denote the set of those measures $\hat{\nu}\in\mathcal M_1(\mathcal M_1(X))$ that $f$ statistically bifurcates toward $\hat{\nu}$ through perturbations in $\Lambda$.
\begin{definition}
    A map $f\in\Lambda$ is called statistically unstable in law w.r.t $\Lambda$ if the set $B_{\Lambda,f}$ has more than one element. 
\end{definition}
\begin{remark}\label{Rem.acc.subset}
By definition, it holds true that 
$$acc(\{\hat{e}_n^f\}_n)\subset B_{\Lambda,f}.$$
\end{remark}
Here are some nice properties of the set $\mathcal{B}_{\Lambda,f}$:
\begin{lemma}\label{compact BLambda}
The set $\mathcal{B}_{\Lambda,f}$ is a compact subset of $\mathcal M_1(\mathcal M_1(X))$. 
\end{lemma}
\begin{proof}
By the definition, it is clear that the set $\mathcal{B}_{\Lambda,f}$ is closed. The compactness is a consequence of the compactness of $\mathcal M_1(\mathcal M_1(X))$.
\end{proof}
\begin{lemma}
The set $\mathcal{B}_{\Lambda,f}$ is a connected subset of $\mathcal M_1(\mathcal M_1(X))$.
\end{lemma}
\begin{proof}
For the sake of the contrary assume that $B_{\Lambda,f}$ is not connected, and can be decomposed to two non-empty disjoint closed sets $A$ and $B$. Therefore there is some real number $d>0$ such that $\tilde{d}_w(A,B)>d.$ Take two elements $\hat{\nu}\in A$ and $\hat{\eta}\in B$ and let $N$ in Lemma \ref{lemma.dec.dist} is chosen so that 
$$\forall n>N, \hat{d}_w(\hat{e}_n^f,\hat{e}_{n+1}^f)<\frac{d}{3}.$$
We can find a neighborhood $U$ of $f$ so that 
$$\forall g\in U, \quad \hat{d}_w(\hat{e}_n^f,\hat{e}_{n}^g)<\frac{d}{3}.$$
This is possible since the map sending $f$ to $\hat{e}_N^f $ is continuous. Now take two maps $h,g\in U$ such that for some integers $n_1,n_2>N$ it holds true that
$$\hat{d}_w(\hat{e}_{n_1}^g,\hat{\nu})<\frac{d}{3} \quad \text{and}\quad \hat{d}_w(\hat{e}_{n_2}^h ,\hat{\eta})<\frac{d}{3}.$$

Consider the following sequence of elements of $\mathcal M_1(\mathcal M_1(X))$: 
$$\hat{\nu}, \hat{e}_{n_1}^g,\hat{e}_{n_1 -1}^g,...,\hat{e}_{N}^g,\hat{e}_{N}^h,...,\hat{e}_{n_2 -1}^h,\hat{e}_{n_2}^h,\hat{\eta}.$$
The distance between two consecutive elements of this sequence is less than $\frac{d}{3}$, and hence there is an element of this sequence that lies out of $\frac{d}{3}$ neighborhood of $A\bigcup B=B_{\Lambda,f}$. By taking $N$ larger, we obtain another element of
 $\mathcal M_1(\mathcal M_1(X))$ out of $\frac{d}{3}$ neighborhood of $B_{\Lambda,f}$. So there is a sequence like $\hat{e}_{n_k}^{f_k}$ out of $\frac{d}{3}$ neighborhood of $B_{\Lambda,f}$ and because of the compactness of $\mathcal M_1(\mathcal M_1(X))$ this sequence has a convergent subsequence converging to an element out of $B_{\Lambda,f}$. By definition, any accumulation point of this sequence is an element of $B_{\Lambda,f}$ which is a contradiction.
 \end{proof}
 
\begin{lemma}
For any $\hat{\nu }\in \mathcal B_{\Lambda,f}$, any measure $\nu$ in the support of $\hat{\nu}$ is invariant under iteration of $f$.
\end{lemma}
\begin{proof}
By definition there is a sequence of maps $\{f_k\}_k$ in $\Lambda$ converging to $f$ and a sequence of natural numbers $\{n_k\}_k$ converging to infinity such that 
$$\lim_{k\to\infty} \hat{d}_w(\hat{e}_{n_k}(f_k),\hat{\nu})=0.$$
If $\nu$ is in the support of $ \hat{\nu}$ then for any neighbourhood $\mathcal U$ of $\nu$ and for $k$ large enough, we have 
$$\hat{e}_{n_k}(f_k)(\mathcal U)>0.$$
Recalling that 
$$\hat{e}_{n_k}(f_k)(\mathcal U)=\int_{X}\delta_{e^{f_k}_{n_k}(x)}(\mathcal U)d\mu,$$
we conclude that the integrand of the integral above is non-zero for a subset of $X$ with positive measure and hence in particular for each $k$ there is a point $x_k\in X$ such that $e^{f_k}_{n_k}(x_k)\in \mathcal U$. Since $\mathcal U $ is an arbitrary neighbourhood of $\nu$ we can choose $x_k$ such that 
$$\lim_{k\to\infty} e_{n_k}^{f_k}(x_k)=\nu.$$
On the other hand note that for a large $k$ the map $f_k$ is close to the map $f$ so the measure $e^{f_k}_{n_k}(x_k)$ is close to $f_*(e^{f_k}_{n_k}(x_k)$. So we have 
$$\lim_{k\to\infty}d_w(e^{f_k}_{n_k}(x_k),f_*(e^{f_k}_{n_k}(x_k)))=0,$$
which together with the continuity of $f_*$ imply that $f_*(\nu)=\nu$ and so we are done.
\end{proof}
The set $\mathcal B_{\Lambda,f}$ depends on the choice of the set of dynamics $\Lambda$ in which we are allowed to perturb the map $f$. If $\Lambda $ is replaced by a larger set of maps, then we may have more elements in $\mathcal M_1(\mathcal M_1(X))$ toward which $f$ statistically bifurcates. 
In this Appendix we will show that:
\begin{customthm}{E}\label{thm.generic BLambda}
 Suppose $X$ is a compact metric space endowed with a reference probability measure $\mu$, and $\Lambda$ is a set of continuous self-mappings of $X$, endowed with a topology finer than $C^0$ topology. Then for Baire generic maps $f\in\Lambda$, we have $$B_{\Lambda,f}=acc(\{\hat{e}_n^f\}_n).$$
 \end{customthm}

 \begin{remark}
     Theorem \ref{thm.generic BLambda} is the counterpart of Theorem \ref{thm.main.essential} in the version of convergence in law.
 \end{remark}
We will prove this theorem later. Let us remind some definitions that we need in the rest of this section. Let $X$ and $Y$ be two topological spaces with $Y$ compact. Denote the set of all compact subsets of $Y$ by $C(Y)$.
\begin{definition}
 A map $\phi:X\to C(Y)$ is called lower semi-continuous if for any $x\in X$ and any $V$ open subset of $Y$ with $\phi(x)\cap V\neq \emptyset$, there is a neighborhood $U$ of $x$ such that for any $y\in U$ the intersection $\phi(y)\cap V$ is non-empty. The map $\phi$ is called upper semi-continuous if, for any $x\in X$ and any $V$ open subset of $Y$ with $\phi(x)\subset V$, there is a neighborhood $U$ of $x$ such that for any $y\in U$ the set $\phi(y)$ is contained in $V$. And finally, $\phi$ is called continuous at $x$ if it is both upper and lower semi-continuous at $x$.
\end{definition}

We also recall the following theorem of Fort \cite{fort1951points} which generalizes the well-known theorem about real-valued semi-continuous maps to the set-valued semi-continuous maps:
\begin{theorem*}[Fort]
For any Baire topological space $X$ and compact topological space $Y$, the set of continuity points of a semi-continuous map from $X$ to $C(Y)$ is a Baire generic subset of $X$.
\end{theorem*}
We recall the following fact on the semi-continuity property of the map sending the dynamics to its set of invariant probability measures. 

\begin{lemma}\label{semi-cont inv meas}
The map sending $f\in \Lambda$ to its set of invariant probability measures is upper semi-continuous. 
\end{lemma}

We recall that by Lemma \ref{compact BLambda}, the set $\mathcal B_{\Lambda,f}$ is compact. We can ask about the dependence of the set $\mathcal B_{\Lambda,f}$ on the map $f$. The following lemma shows that this dependency is semi-continuous:
\begin{lemma}\label{lemma.BLambda.usc}
The map sending $f\in\Lambda $ to the set $\mathcal B_{\Lambda, f}$ is upper semi-continuous.
\end{lemma}
For the proof see Lemma 1.10. in \cite{talebi20}.\\

To each map, $f\in\Lambda$ one can associate the set of accumulation points of the sequence $\{\hat{e}_n(f)\}_{n\in \mathbb{N}}$ which is a compact subset of $\mathcal M_1 (\mathcal M _1(X))$. This map is neither upper semi-continuous nor lower semi-continuous. However if we add the points of this sequence to its accumulation points and consider the map sending $f\in \Lambda$ to the closure $\overline{\{\hat{e}_n(f)|n\in \mathbb{N}\}}$, we obtain a semi-continuous map:
\begin{lemma}\label{E lsc}
The map $\mathcal E: \Lambda \to C({\mathcal M_1(\mathcal M_1(X))})$ defined as 
$$\mathcal E(f):=\overline{\{\hat{e}_n(f)|n\in \mathbb{N}\}},$$
is lower semi-continuous.
\end{lemma}
For the proof see Lemma 1.11. in \cite{talebi20}.\\

The following lemma is an interesting consequence of Lemma \ref{E lsc} that shows how the set $\mathcal{E}(f)$ depends on the dynamics $f$. 

\begin{lemma}\label{cor of Fort}
A Baire generic map $f\in\Lambda$ is a continuity point for the map $\mathcal E$.
\end{lemma}

This lemma gives a view to the statistical behaviors of generic maps in any Baire space of dynamics: for a generic map, the statistical behavior that can be observed for times close to infinity can not be changed dramatically by small perturbations.  
\begin{proof}
Using Lemma \ref{E lsc}, this is a direct consequence of Fort's theorem. 
\end{proof}
Now we are ready to prove Theorem \ref{thm.generic BLambda}. This theorem reveals how two notions of statistical instability in law and non-statistical maps in law are connected. There is another proof of this theorem which is communicated by Pierre Berger that can be found in \cite{talebi20} (Theorem 1.14).

 \begin{proof}[Proof of Theorem \ref{thm.generic BLambda}]
Take a generic map $f$ from the main lemma above. By remark \ref{Rem.acc.subset}, $acc(\{\hat{e}_n^f\}_n)\subset B_{\Lambda,f}$. So if $B_{\Lambda,f}$ has only one element, since $acc(\{\hat{e}_n^f\}_n)$ is non-empty, the equality holds. Now suppose that $B_{\Lambda,f}$ has more than one element. For the sake of contradiction suppose there is a measure $\hat{\nu}\in B_{\Lambda,f}$ which is not in $acc(\{\hat{e}_n^f\}_n)$. Then there is a number $n\in\mathbb{N}$ such that $\hat{\nu}=\hat{e}_n^f$ and $\hat{e}_n^f$ is an isolated point of the sequence $\{\hat{e}_n^f\}_n$. Recalling that for generic $f$ we have $B_{\Lambda,f}\subset \mathcal{E}(f)$, we can conclude that $B_{\Lambda,f}$ can be written as a union of two disjoint and non-empty closed set:
$$B_{\Lambda,f}= \{\hat{e}^f_n\} \bigcup \{\hat{e}^f_n\}^c.$$
This is contrary to the connectedness of the set $B_{\Lambda,f}$. 
 \end{proof}

Consider all invariant measures $\nu$ that the map $f$ statistically bifurcates toward $\hat{\delta}_nu$. We denote these measures by $\mathcal M_{\Lambda,f}$ which are defined more precisely as bellow:
$$\mathcal M_{\Lambda, f}:=\{\nu\in \mathcal M_1(X)|\hat{\delta}_{\nu}\in\mathcal B_{\Lambda,f}\}.$$
\begin{theorem}\label{generic stat behav}
Let $\Lambda$ be a Baire space of self-mappings of $X$ endowed with a topology finer than $C^0$-topology. For a Baire generic map $f\in\Lambda$ the empirical measures of $\mu$ almost every point $x\in X$, accumulates to each measure in $\mathcal M _{\Lambda,f}$ or in other words:
\begin{align}\label{genericity of large accumulation}
for\ \mu-a.e. \ x\in X,\quad \mathcal M_{\Lambda,f}\subset acc(\{e^f_n(x)\}_{n\in\mathbb N}).
\end{align}
\end{theorem}
For the proof see Theorem 1.16. in \cite{talebi20}.\\

If one can find any information about the set $\mathcal M_{\Lambda,f}$ for a generic map $f$ in $\Lambda$ then by Theorem \ref{generic stat behav}, we can translate this information to information about the statistical behavior of $\mu$-almost every point for a generic subset of maps. 

The following lemma shows how the set $\mathcal M_{\Lambda,f}$ depends on the map $f$:
\begin{lemma}\label{MLambda,f usc}
The map sending $f\in\Lambda $ to the set $\mathcal M_{\Lambda, f}$ is upper semi-continuous.
\end{lemma}
For the proof see Lemma 1.16. in \cite{talebi20}.\\

Now let us see what is the consequence of this lemma and Theorem \ref{generic stat behav} when we know the maps in a dense subset bifurcates toward each Dirac mass at an invariant measure.

\begin{proposition}[Maximal oscillation]\label{maximal-oscil}
If there is $D\subset\Lambda$ dense such that any map $f\in D$ bifurcates toward the Dirac mass at each invariant measure through perturbations in $\Lambda$, or in another word $M_{\Lambda,f}=\mathcal M_1(f)$, then a generic $f\in\Lambda$ has maximal oscillation. 
\end{proposition}
For the proof see Proposition 1.19. in \cite{talebi20}.\\

\section{Comparison between different versions}\label{sec.comparison}

 In this section, we compare different versions of defining statistical instability and non-statistical dynamics and show how they are related. 
 The first proposition describes the relation between different versions of defining non-statistical maps.
 \begin{proposition}\label{prop hier non-stat in law}
 Suppose $f:X\to X$ is a continuous map of a compact metric space. \\
 i) If $f$ is non-statistical in law, then it is $L^1$ non-statistical. \\
 ii) If $f$ is $L^1$ non-statistical, then it is non-statistical. 
  \end{proposition}
  
  \begin{proof}
 i) Let $f$ be non-statistical in law. So by definition the sequence $\{\hat{e}_n(f)=(e^f_n)_* \mu\}_n$ is not convergent. We recall that 
 \[(e_n^f)_* (\mu)=\int_X\hat{\delta}_{e^f_n(x)}d\mu(x).\]
 where $\hat{\delta}_{e^f_n(x)}\in \mathcal{M}_1(\mathcal{M}_1(X))$ is the Dirac mass supported on the point $e^f_n(x)\in \mathcal{M}_1(X)$.\\
 Suppose to the contrary that $f$ is not $L^1$ non-statistical. So the sequence of maps $\{e_n^f:X\to \mathcal{M}_1(X)\}$
  is convergent in the $L^1$ topology. Let us call the limit point of this sequence by $e^f_\infty$. Now we show that $(e^f_n)_*\mu$ converges to $(e^f_\infty)_*\mu$ which is a contradiction. For simplicity we denote $(e_n^f)_*\mu$ by $\nu_n$ and $(e_\infty^f)_*\mu$ by $\nu_\infty$. We recall that 
  \[\hat{d}_w(\nu_n,\nu_\infty)=\min\limits_{\xi\in\pi(\nu_n,\nu_\infty)}\int_{\mathcal{M}_1(X)\times \mathcal{M}_1(X)} d_w(e,e') d\xi(e,e').\]
  where $\pi(\nu_n,\nu_\infty)$ is the set of all probability measures on $\mathcal{M}_1(X)\times \mathcal{M}_1(X)$ which projects to $\nu_n$ and $\nu_\infty$ under the projections to the first and second coordinate respectively. Consider the following element of $\pi(\nu_n,\nu_\infty)$:
  \[\xi:=\int_Xd_w(\delta_{e_n^f(x)},\delta_{e_\infty^f(x)})d\mu\]
  We have 
   $$\hat{d}_w(\nu_n,\nu_\infty)\leq \int_{\mathcal{M}_1(X)\times \mathcal{M}_1(X)}d_w(e,e')d\xi_n(e,e').$$
  On the other hand, 
   \[\int_{\mathcal{M}_1(X)\times \mathcal{M}_1(X)}d_w(e,e')d\xi_n(e,e')=\int_Xd_w(e^f_n(x), e^f_\infty(x))d\mu=d_{L'}(e_n^f, e_\infty^f).\]
  So we obtain $\hat{d}_w(\nu_n,\nu_\infty)\leq d_{L^1}(e_n^f,e_\infty^f)$, which implies that $\nu_n$ is converging to $\nu_\infty$ and is a contradiction. \\
   
   ii) Let $f$ be $L^1$ non-Statistical. We show that $f$ is non-statistical and so the maps $e_n^f:X\to \mathcal{M}_1(X)$ do not converge almost surely. By contrary, suppose the maps $e_n^f:X\to \mathcal{M}_1(X)$ converge almost surely to a map $e_\infty^f:X\to \mathcal{M}_1(X)$. Hence the map $d_w(e_n^f(.), e_\infty^f(.)):X\to \mathbb{R}$ converges to zero almost surely, and by dominated convergence theorem we obtain that 
   \[ d_{L^1}(e_n^f,e_\infty^f)=\int_X d_w(e_n^f(x),e^f_\infty(x))d\mu(x)\to 0 \quad (n\to \infty).\]
   Which is a contradiction. 
  \end{proof}

  We would like to announce here the existence of examples of maps that are $L^1$ statistical but non-statistical and the maps that are statistical in law but $L^1$ non-statistical. These examples could be constructed on $\mathbb S^1\times \mathbb A$, where $\mathbb A$ is the closed annulus, using the Anosov-Katok method and similar ideas described in Section \ref{sec.Anosov-Katok}. \\
  
  The next proposition shows that the same hierarchy holds for different versions of defining statistically unstable maps.
  \begin{proposition}\label{prop hier stat stab}
  Suppose $\Lambda$ is a set of continuous self-mappings of a compact metric space $X$. It holds true that:\\
  i) If $f\in \Lambda$ is statistically unstable in law, then it is $L^1$ statistically unstable. \\
  ii) If $f\in \Lambda$ is $L^1$ statistically unstable, then it is statistically unstable. 
  \end{proposition}
    \begin{proof}
    i) If $f$ is statistically in law, then by definition there are at least two different elements $\hat{\nu}_1$ and $\hat{\nu}_2$ in $\mathcal{M}_1(\mathcal{M}_1(X))$ toward which $f$ statistically bifurcates. This means that there are two sequences of maps $\{f^1_k\}_k$ and $\{f^2_k\}_k$ converging to $f$, and two sequences of positive integers $\{n^1_k\}_k$ and $\{n^2_k\}_k$ converging to infinity such that 
    \[\lim\limits_{k\to \infty} \hat{e}_{n^1_k}(f^1_k) =\hat{\nu}_1, \lim\limits_{k\to \infty} \hat{e}_{n^2_k}(f^2_k) =\hat{\nu}_2\]
    Now on the contrary suppose $f$ is $L^1$ statistically stable. So the sequence $\{e^{f^1_k}_{n^1_k}:X\to X\}_k$ and $\{e^{f^2_k}_{n^2_k}:X\to X\}_k$ 
    both converge to a map $e^f_\infty:X\to X$ in the $L^1$ topology. Using the same arguments in the proof of part (i) in the previous proposition we conclude that both sequences $\{(e^{f^1_k}_{n^1_k})_* \mu\}_k$ and $\{(e^{f^2_k}_{n^2_k})_* \mu\}_k$
    converge to $(e_\infty^f)_*\mu$. Hence we have $\hat{\nu}_1=\hat{\nu}_2=(e_\infty^f)_*\mu$ and since $\hat{\nu}_1$ and $\hat{\nu}_2$ were distinct elements $\mathcal{M}_1(\mathcal{M}_1(X))$, this is a contradiction.\\
    ii) Suppose $f$ is not statistically stable, so there is a map $e_\infty^f:X\to \mathcal{M}_1(X)$ such that for any sequence $\{f_k\}_k$ converging to $f$ and any sequence of natural numbers $\{n_k\}_k$ converging to infinity the sequence of maps $\{e_{n_k}^{f_k}:X\to \mathcal{M}_1(X)\}_k$ converge almost surely to the map $e_\infty^f$. Using the dominated convergence theorem, we conclude the convergence of this sequence in the $L^1$ topology to the map $e_\infty^f$ and hence $f$ can not be $L^1$ statistically unstable. 
        \end{proof}

\def\polhk#1{\setbox0=\hbox{#1}{\ooalign{\hidewidth
  \lower1.5ex\hbox{`}\hidewidth\crcr\unhbox0}}}

\bibliographystyle{plain}
\bibliography{References}

\end{otherlanguage}
\end{document}